\documentclass[12pt]{amsart}

\usepackage[active]{srcltx}
\usepackage{graphics,latexsym,dsfont,pictex,dcpic}
\usepackage{latexsym}
\usepackage{amsmath,amsthm,amssymb,bbm}
\usepackage{mathrsfs}

\newcommand{\ca}[1]{\mathcal{#1}}

\newcommand{\bb}[1]{\mathbb{#1}}
\newcommand{\bbm}[1]{\mathbbm{#1}}
\newcommand{\fr}[1]{\mathfrak{#1}}
\renewcommand{\sc}[1]{\mathscr{#1}}

\newcommand{\Naturals}{\mathds{N}}
\newcommand{\Integers}{\mathds{Z}}
\newcommand{\Rationals}{\mathds{Q}}
\newcommand{\Reals}{\mathds{R}}
\newcommand{\Complex}{\mathds{C}}

\newcommand{\id}{\mathrm{id}}
\newcommand{\Symp}{\mathrm{Sp}}
\newcommand{\Unitary}{\mathrm{U}}
\newcommand{\Lie}{\mathrm{Lie}}
\newcommand{\GeneralLinear}{\mathrm{GL}}
\newcommand{\SpecialUnitary}{\mathrm{SU}}
\newcommand{\SpecialLinear}{\mathrm{SL}}
\newcommand{\SpecialLinearLie}{\mathfrak{sl}}
\newcommand{\SpecialOrthogonalLie}{\mathfrak{so}}
\newcommand{\SpecialUnitaryLie}{\mathfrak{su}}
\newcommand{\SymplecticGroup}{\mathrm{Sp}}
\newcommand{\SymplecticLie}{\mathfrak{sp}}
\newcommand{\PSpc}{\mathbb{P}}
\newcommand{\ASpc}{\mathbb{A}}
\newcommand{\Proj}{\mathrm{Proj}}
\newcommand{\ModularCY}{\fr{M}_{\mathrm{CY}}}
\newcommand{\HyperCY}{\fr{H}_{\mathrm{CY}}}
\newcommand{\Cohom}{\mathrm{H}}
\newcommand{\Jacobian}{\mathrm{Jac}}
\newcommand{\Hom}{\mathrm{Hom}}

\newcommand{\Aut}{\mathrm{Aut}}

\newcommand{\IdMatrix}{\mathbbm{1}}
\newcommand{\NormalForm}{\mathrm{NF}}
\newcommand{\transpose}[1]{{}^{\mathrm{t}} #1}

\theoremstyle{plain}
\newtheorem{thm}{Theorem}[section]
\newtheorem{theorem}[thm]{Theorem}
\newtheorem{lemma}[thm]{Lemma}

\newtheorem{proposition}[thm]{Proposition}

\theoremstyle{definition}

\newtheorem{remark}[thm]{Remark}

\numberwithin{equation}{thm}

\title[Disproof of Modularity]{
    Disproof of modularity of moduli space of CY 3-folds of double covers of
     $\mathbb P^3$ ramified along eight planes
     in general positions}

\author[Ralf Gerkmann]{Ralf Gerkmann}
\address{Universit\"at Mainz, Fachbereich 17, Mathematik,
    55099 Mainz, Germany}
\email{gerkmann@mathematik.uni-mainz.de}

\author[Mao Sheng]{Mao Sheng${}^\dagger$}
\address{East China Normal University, Dep.\ of Mathematics, 200062
    Shanghai, P.R.\ China}
\email{msheng@math.ecnu.edu.cn}
\thanks{${}^\dagger$ The second named author is supported by a postdoctoral fellowship in the East China Normal University.}

\author[Kang Zuo]{Kang Zuo}
\address{Universit\"{a}t Mainz, Fachbereich 17, Mathematik, 55099 Mainz,
    Germany}
\email{kzuo@mathematik.uni-mainz.de}

\begin{document}

\fontsize{10}{14}\selectfont \maketitle

\begin{abstract}
We prove that the moduli space of Calabi-Yau 3-folds coming from
eight planes of $\PSpc^3$ in general positions is not modular. In
fact we show the stronger statement that the Zariski closure of the
monodromy group is actually the whole $\SymplecticGroup(20,\Reals)$. We
construct an interesting submoduli, which we call
\emph{hyperelliptic locus}, over which the weight 3
$\Rationals$-Hodge structure is the third wedge product of the
weight 1 $\Rationals$-Hodge structure on the corresponding
hyperelliptic curve. The non-extendibility of the hyperelliptic locus
inside the moduli space of a genuine Shimura subvariety is proved.
\end{abstract}

\section{Introduction}
\label{sec:Intro}

In the study of geometry of moduli space, it is important to
characterize those moduli spaces which are locally Hermitian
symmetric varieties. We refer the reader to \cite{VZ1},
\cite{VZ2}, \cite{MVZ} for such a theory based on the Arakelov
equality. On the other hand, in order to prove a negative result
it is also important to find some necessary conditions, which can
be checked quite easily for explicitly given moduli spaces. In
this paper, we will work with an interesting moduli space of CY
3-folds, which comes from the hyperplane arrangements in $\PSpc^3$
consisting of eight planes in general positions. The aim of our
present work is to disprove the modularity of this moduli space by
two different methods. Before stating the main theorem, we shall
make the meaning of modularity precise, since it could be
ambiguous in certain cases. For example, the moduli space of six
lines of $\PSpc^2$ in general positions, which is identical to the
moduli space of six points of $\PSpc^2$ in general positions, can
be openly embedded either into an arithmetic quotient of type four
bounded symmetric domain \cite{MSY} or into an arithmetic ball
quotient \cite{ACT}, \cite{DKV} by different period mappings. \\

Let $\fr{M}$ be the coarse moduli scheme representing a moduli
functor $\ca{M}$ of polarized algebraic manifolds of dimension
$n$. After a finite base change of $\fr{M}$, one obtains a
universal family $f: \ca{X} \to S$. The rational primitive middle
cohomologies of the fibers of $f$ constitute a
$\Rationals$-polarized variation of Hodge structures $\bb{V}$. It
induces the period mapping
$$
\phi: S\hookrightarrow \Gamma\backslash G'/K'
$$
where $G'/K'$ is the classfying space of polarized Hodge
structures. If there exists a locally Hermitian symmetric variety
$\Gamma\backslash G/K$ and a locally homogenous PVHS $\bb{W}$ over
it such that $\phi$ factors through the period mapping
$$
\psi: \Gamma\backslash G/K \hookrightarrow \Gamma'\backslash G'/K'
$$
defined by $\bb{W}$ and such that the induced map
$$
\phi: S \to \Gamma\backslash G/K
$$
is an open embedding, then we say that $f$ is a {\em modular
family} with respect to $(G/K,\bb{W})$. In the case that the
reference locally homogenous PVHS $\bb{W}$ and the Hermitian
symmetric space $G/K$ are clear from the context, we simply say
$f$ is modular family. The moduli space $\fr{M}$ is said to be
{\em modular} if a certain universal family of $\fr{M}$ is a
modular family. Under this definition, it is clear that if
$$
f: \ca{X} \rightarrow S
$$
is a modular family with respect to $(G/K,\bb{W})$, then one has a
factorization of the monodromy representation $\rho$ of $\bb{V}$:

\vspace*{0.5cm}\hspace*{3.5cm}
\begindc{\commdiag}[25]
\obj(1,3){$\pi_1(S)$} \obj(5,3){$G'$} \obj(3,2){$G$}
\mor{$\pi_1(S)$}{$G'$}{$\scriptstyle\rho$} \mor{$\pi_1(S)$}{$G$}{}
\mor{$G$}{$G'$}{$\scriptstyle\varrho$}[-1,0]
\enddc
\vspace*{0.3cm}

\noindent where $\varrho: G\to G'$ is the group homomorphism
determined by
$\bb{W}$. \\

Now let $\ModularCY$ be the moduli space of CY 3-folds from eight
planes of $\PSpc^3$ in general positions. The classfying space of
the polarized Hodge structure on the middle cohomology of such a
CY 3-fold is
$$
D'=\frac{\Symp(20,\Reals)}{\Unitary(1)\times \Unitary(9)}.
$$
The natural Hermitian symmetric space in this case is
$$
D=\frac{\SpecialUnitary(3,3)}{\mathrm{S}(\Unitary(3)\times
\Unitary(3))} \quad,
$$
and the locally homogenous PVHS $\bb{W}$ is the Calabi-Yau like
PVHS over $\Gamma\backslash D$ \mbox{(cf.\ \cite{SZ})}, which is
induced from the group homomorphism
$$
\bigwedge^{3}: \SpecialUnitary(3,3)\to \Symp(20,\Reals).
$$
The main theorem of this paper is the following

\begin{theorem}\label{MainTheorems}
Let $\ModularCY$ and $(D,\bb{W})$ be as above. Then $\ModularCY$
is not modular with respect to $(D,\bb{W})$.
\end{theorem}

To keep our theorem in perspective, we would like to point out
that the analogous moduli spaces of CY $n$-folds for $n\leq 2$ are
modular (cf.\ \cite{MSY}). It would be very interesting to
extend the present work to the $n\geq 4$ cases. At this point, we
would like to remind the reader of the early work \cite{SYY}. They
disproved a modularity result similar to Theorem
\ref{MainTheorems} in a more general setting. It seems that the
theory used in \cite{SYY} has not been widely accepted within the
community of algebraic geometricians. We hope our purely Hodge
theoretical proof will at least clarify some serious issues about
the disproof of modularity. Furthermore, by a study of possible real
groups of Hodge type contained in $\Symp(20,\Reals)$ and an
application of the plethysm method we can deduce a stronger result
about the Zariski closure of the monodromy group.

\begin{theorem}\label{ZarClosMon}
Let $f: \ca{X} \rightarrow S$ be a universal family of
$\ModularCY$, and $\bb{V}=R^{3}f_{*}\Rationals_\ca{X}$ be the
interesting weight 3 VHS. Then the image of the monodromy
representation $\rho: \pi_{1}(S)\to \Symp(20,\Reals)$ of $\bb{V}$
is Zariski dense.
\end{theorem}

In the work of \cite{MSY}, a special submoduli, which is
isomorphic to an arithmetic quotient of
$\Symp(4,\Reals)/\Unitary(2)$, was constructed. We shall
generalize their construction to our case. Since this submoduli
arises from the moduli of hyperelliptic curves of genus 3, we simply
call it the \emph{hyperelliptic locus}. We show that over the
five dimensional hyperelliptic locus the weight 3 VHS of CY
3-folds is isomorphic to the wedge product of weight 1 VHS (cf.\
Prop.\ \ref{wedge product of weight one HS}). It is then natural to
ask if one can extend the hyperelliptic locus in $\ModularCY$ to a
six dimensional submoduli which is isomorphic to an arithmetic
quotient of $\Symp(6,\Reals)/\Unitary(3)$. Using the concept of
characteristic subvarieties we arrive at a negative answer of this
extension problem.

\begin{theorem}\label{NonExtHyperLocus}
Let $\HyperCY$ be the hyperelliptic locus of $\ModularCY$. Then
there exists no extension of $\HyperCY$ inside $\ModularCY$ such
that it is isomorphic to a Zariski open subset of an arithmetic
quotient of $\Symp(6,\Reals)/\Unitary(3)$ which is embedded into
$\ModularCY$ via the wedge product of the weight 1 VHS of a
universal family of abelian 3-folds.
\end{theorem}

The paper is organized as follows. In \S\ref{sec:CYConstruction}
we will construct two different Calabi-Yau manifolds from a given
hyperplane arrangement, and describe the relation between them.
The construction of the hyperelliptic locus concludes the second
section. In \S\ref{sec:MethodDescr} we will describe our
methodology to disprove the modularity. Two different methods will
be presented respectively. The actual computations for our moduli
space are realized using the theory of Jacobian rings. We have to
adapt the current knowledge of Jacobian ring to our case. This is
done in the fourth section. Section \ref{sec:Proofs} contains the
results of our computations and the proof of the main theorems
stated in \S\ref{sec:Intro}. \\[0.2cm]

\noindent {\em Acknowledgements:} The major part of this work was
done during an academic visit of the second named author to the
department of mathematics of the university of Mainz in 2006. He
would like to express his hearty thanks to the hospitality of the
faculty, especially Stefan M\"{u}ller-Stach. We also thank Eckart
Viehweg for his interests and several helpful conversations about
this work.

\section{Calabi-Yau manifolds from eight planes of $\PSpc^3$ in general
    positions}
\label{sec:CYConstruction}

Let $H_1,\cdots,H_8$ denote eight planes of $\PSpc^3$ in general
position. They sum up to a simple normal crossing divisor
$$
B=H_1+\cdots+H_8
$$
on $\PSpc^3$. Since $B$ is even, we can form the double cover $X$
of $\PSpc^3$ with branch locus $B$. Obviously, $X$ is only a
singular variety because $B$ is singular. Fixing an ordering of
the irreducible components $H_{ij}=H_{i}\cap H_{j}$ of
singularities of $B$, we use the canonical resolution of double
covers to obtain a smooth model $\tilde{X}$ of $X$. Namely, we
have the following commutative diagram

\vspace*{0.5cm}\hspace*{5.0cm}
\begindc{\commdiag}[25]
\obj(1,3){$X$} \obj(3,3)[Xs]{$\tilde{X}$} \obj(1,1)[P3]{$\PSpc^3$}
\obj(3,1)[Ps]{$\tilde{\PSpc}^3$}
\mor{Xs}{$X$}{$\scriptstyle\tau$}[-1,0]
\mor{Ps}{P3}{$\scriptstyle\sigma$}[-1,0]
\mor{$X$}{P3}{$\scriptstyle\pi$}[-1,0]
\mor{Xs}{Ps}{$\scriptstyle\tilde{\pi}$}
\enddc
\vspace*{0.3cm}

\noindent where $\sigma: \tilde{\PSpc}^3 \rightarrow \PSpc^3$ is
the composition of the sequence of blows-up with smooth centers
(the strict transform of) $H_{ij}$. The variety $\tilde{X}$ is a
smooth projective CY 3-fold with
$$
h^{2,1}(\tilde{X})=9 \qquad \text{and} \qquad
h^{1,1}(\tilde{X})=29.
$$
This construction can actually be extended to all $2n+2$
hyperplane arrangements of $\PSpc^n$ in general positions, and the
Hodge numbers of the primitive middle cohomology of the resulting
smooth CY $n$-fold $\tilde{X}$ are
$$
h^{p,n-p}_{pr}(\tilde{X})=\binom{n}{p}^{2}.
$$
For the details, we refer to chapter 3 of \cite{S}.\\

It is easy to see that the moduli space of ordered eight
hyperplane arrangements of $\PSpc^3$ in general positions is of
dimension 9. Hence, by fixing an ordering of the index set
$$
I=\{ (i,j) \in \Naturals^2 ~\mid~ 1\leq i<j\leq 8\} ~~,
$$
the above constructions gives rise to a complete moduli scheme
$\ModularCY$ of smooth CY 3-folds. We note that a different
ordering of $I$ yields a different birational minimal model of the
singular CY $X$.\\

Now we consider the embedding determined by the starting
hyperplane arrangement $(H_1,...,H_8)$, namely
$$
j : \PSpc^3 \hookrightarrow \PSpc^7, \qquad x \mapsto
(\ell_{1}(x): ~\cdots~ :\ell_{8}(x))
$$
where $\ell_i : \Complex^4 \rightarrow \Complex$ denotes a linear
form such that $H_{i}=\{ x \in \PSpc^3 \mid \ell_i=0 \}$ for $1
\leq i\leq 8$. The defining equations of $j(\PSpc^{3})\subset
\PSpc^7$ are the four linearly independent relations among
$\ell_{1},\cdots,\ell_{8}$. Written out explicitly, they are
$$
\begin{array}{ccccccc}
a_{11}y_1 & + & \cdots & + & a_{18}y_8 & = & 0 \\
 \vdots & & & & \vdots \\
a_{41}y_1 & + & \cdots & + & a_{48}y_8 & = & 0
\end{array}
$$
where $(y_1:~\cdots~:y_8)$ denote homogenous coordinates of
$\PSpc^7$. We can define a new CY 3-fold $Y$ which is the complete
intersection of four quadrics in $\PSpc^7$ defined by
\begin{equation}\label{FourQuadrics}
\begin{array}{ccccccc}
a_{11}y_1^2 & + & \cdots & + & a_{18}y_8^2 & = & 0 \\
 \vdots & & & & \vdots \\
a_{41}y_1^2 & + & \cdots & + & a_{48}y_8^2 & = & 0.
\end{array}
\end{equation}
The variety $Y$ is smooth since any $4\times 4$ minor of the
matrix $A := (a_{ij})$ is nonzero. The two Hodge numbers of $Y$
are computed to be
$$
h^{2,1}(Y)=65 \qquad \text{and} \qquad h^{1,1}(Y)=1.
$$
In particular, the moduli space of complex structures on $Y$ is of
dimension $65$. The covering map
$$
\PSpc^7 \longrightarrow \PSpc^7,\qquad (y_1:~\cdots~:y_8) \mapsto
(y_1^2:~\cdots~:y_8^2)
$$
restricts to $p: Y \rightarrow j(\PSpc^3)$. The composite map,
denoted again by $p$,
$$
p : Y \longrightarrow j(\PSpc^3) \cong \PSpc^3
$$
exhibits $Y$ as the Kummer covering of $\PSpc^3$, branched along
$B$ with degree $2^7$. Clearly, the Galois group $\Aut(Y|\PSpc^3)$
of $Y$ over $\PSpc^3$ is isomorphic to $(\Integers/2\Integers)^7$.
There is a canonical surjection
$$
G_1 := (\Integers/2\Integers)^8 \longrightarrow \Aut(Y|\PSpc^3),
\qquad a = (a_1,...,a_8) \mapsto \sigma_a
$$
with $\sigma_a(y_i) = (-1)^{a_i} y_i$ for $1\leq i \leq 8$. Its
kernel of order $2$ is generated by $(1,...,1)$. Furthermore, we
have a distinguished index 2 normal subgroup $N_{1} \triangleleft
G_{1}$ given by
$$
N_1 := \ker \left( G_{1} \cong (\Integers/2\Integers)^8
\stackrel{\sum}{\longrightarrow} \Integers/2\Integers \right).
$$
The following proposition reveals the geometric relation between
two CY manifolds coming from the same hyperplane arrangement.

\begin{proposition}\label{XYRelation}
Let $(H_1,...,H_8)$ be a hyperplane arrangement of $\PSpc^3$ in
general position. Then one can construct two smooth CY 3-folds
$\tilde{X}$ and $Y$ as above. One has a natural isomorphism
$$
\Cohom^3(\tilde{X},\Rationals) \cong \Cohom^3(Y,\Rationals)^{N_1}.
$$
\proof The quotient map $p$ factors as
$$
Y \stackrel{\alpha}{\longrightarrow} Y/N_{1} \longrightarrow
\PSpc^3.
$$
The degree of $Y/N_{1}$ over $\PSpc^3$ is $2$ and one can directly
check that $Y/N_{1}$ branches exactly along $B$. Since $\PSpc^3$
has no torsion element in the second integral cohomology, we have
the identification $X=Y/N_1$. Now let us examine the following
commutative diagram

\vspace*{0.5cm}\hspace*{5.0cm}
\begindc{\commdiag}[25]
\obj(1,3){$Y$} \obj(3,3)[Ys]{$\tilde{Y}$} \obj(1,1)[P3]{$\PSpc^3$}
\obj(3,1)[Ps]{$\tilde{\PSpc}^3$}
\mor{Ys}{$Y$}{$\scriptstyle\tilde{\tau}$}[-1,0]
\mor{Ps}{P3}{$\scriptstyle\sigma$}[-1,0]
\mor{$Y$}{P3}{$\scriptstyle p$}[-1,0]
\mor{Ys}{Ps}{$\scriptstyle\tilde{p}$}
\enddc
\vspace*{0.3cm}

\noindent where $\tilde{Y}$ is the normalization of the fiber
product of $Y$ and $\tilde{\PSpc}^3$ over $\PSpc^3$. Obviously,
$\tilde{\tau}$ is a contraction map. Since $Y$ is smooth,
$\tilde{\tau}$ induces the isomorphism
$$
\Cohom^3(Y,\Rationals) \cong \Cohom^3(\tilde{Y},\Rationals).
$$
We put $\tilde{B}$ to be the strict transform of $B$ under
$\sigma$. Then the projection $\tilde{p}$ is the Kummer covering
map of degree $2^7$ with branch locus $\tilde{B}$. Argued as
previously, $\tilde{p}$ factors as
$$
\tilde{Y} \stackrel{\tilde{\alpha}}{\longrightarrow}
\tilde{Y}/N_{1}=\tilde{X}\stackrel{\tilde{p}}{\longrightarrow}
\tilde{\PSpc}^3.
$$
Since $\tilde{\tau}$ is $G_{1}$-equivariant, we have
$\Cohom^3(Y,\Rationals)^{N_{1}} \cong
\Cohom^3(\tilde{Y},\Rationals)^{N_1}$, and as both $\tilde{X}$ and
$\tilde{Y}$ are smooth,
\begin{equation}\label{Hisom}
\Cohom^{3}(\tilde{X},\Rationals)  \cong
\Cohom^3(\tilde{Y},\Rationals)^{N_1}.
\end{equation}
Therefore, combining the last two isomorphisms, we obtain the
isomorphism stated in the proposition.\qed
\end{proposition}

We proceed to construct the hyperelliptic locus $\HyperCY$ inside
our moduli space $\ModularCY$, generalizing the construction in
\cite{MSY}. We first recall that there is a natural Galois
covering
$$
\gamma: (\PSpc^1)^3 \longrightarrow \PSpc^3
$$
with Galois group $S_3$, the symmetric group of three letters.
Explicitly, let
$$
(x_i:y_i), \qquad 1\leq i \leq 3
$$
be the homogenous coordinates of $i$-th factor of $(\PSpc^1)^3$,
such that the components of the quotient map $\gamma$ are given by the
$t$-coefficients of the polynomial
$f(t)=\prod_{i=1}^{3}(x_{i}t+y_{i})$. Now we take arbitrary eight
distinct points $p_1,p_2,...,p_{8} \in \PSpc^1$ and construct a
hyperplane arrangement from it.

\begin{lemma}\label{in general position}
For $1 \leq i \leq 8$ let $H_i$ denote the image of $\{ p_i
\}\times \PSpc^1 \times \PSpc^1$ under $\gamma$. Then $H_i$ is a
hyperplane in $\PSpc^3$, and the hyperplane arrangement
$(H_1,H_2,...,H_8)$ is in general position.

\proof Let $(z_0:z_1:z_2:z_3)$ be the homogenous coordinates of
$\PSpc^3$, and $p=(a:b)$ be a point of $\PSpc^1$. Then using the
expression of $\gamma$, the defining equation of the image set
$\gamma( \{ p \} \times \PSpc^1 \times \PSpc^1 )$ is easily seen
to be
$$
b^{3}z_0-ab^{2}z_1+a^{2}bz_{2}-a^{3}z_3=0.
$$
Therefore, $H_i$ is obviously a hyperplane in $\PSpc^3$. We can
choose an appropriate system of coordinates on $\PSpc^1$ such that
the eight points have coordinates $(-a_i:1)$ for $1\leq i \leq 8$.
Then the columns of the following matrix give the defining
equations of the arrangement $(H_1,H_2,\cdots,H_{8})$:
$$
\left(
  \begin{array}{cccc}
    1 & 1 & \cdots & 1 \\
    a_1&a_2 &\cdots & a_{8} \\
    a_1^2 & a_2^2 &\cdots & a_8^2 \\
    a_1^3 & a_2^3 & \cdots & a_{8}^3 \\
  \end{array}
\right)
$$
Now the property of the hyperplane arrangement to be in general
position is equivalent to that all $4 \times 4$-minors of the
above matrix are nonzero. Since $a_1,a_2,...,a_{8}$ are distinct
from each other by our assumption, all $4 \times 4$-minors are
in Vandermonde form and thus non-zero.
The lemma is proved.\qed
\end{lemma}

Now let $C$ be the hyperelliptic curve over $\PSpc^1$ branched at
$p_1, \cdots,p_{8}$, and let $q$ denote the corresponding covering
map. The Galois group $G_2$ of the composition of morphisms
$$
C^{3} \stackrel{q^3}{\longrightarrow} (\PSpc^1)^3
\stackrel{\gamma}{\longrightarrow} \PSpc^3
$$
is isomorphic to the semi-direct product $N_{2} \rtimes S_3$,
where $N_{2}=\langle\iota_1,\iota_2,\iota_3\rangle$ is the group
generated by the hyperelliptic involutions on each factor of
$C\times C\times C$. One observes that there is a distinguished
index two subgroup $G_{2}'= N_{2}' \rtimes S_3$ of $G_{2}$, where
$N_{2}'$ is the kernel
$$
N_2' := \ker \left( N_{2}\simeq
(\Integers/2\Integers)^3\stackrel{\sum}{\longrightarrow}
\Integers/2\Integers \right)
$$
of the multiplication map. This gives the following commutative
diagram of Galois coverings:

\vspace*{0.5cm}\hspace*{5.0cm}
\begindc{\commdiag}[25]
\obj(1,3)[C3]{$C^3$} \obj(3,3){$X$} \obj(1,1)[P13]{$(\PSpc^1)^3$}
\obj(3,1)[P3]{$\PSpc^3$} \mor{C3}{$X$}{$\scriptstyle\delta$}
\mor{C3}{P13}{$\scriptstyle q^3$} \mor{$X$}{P3}{$\scriptstyle\pi$}
\mor{P13}{P3}{$\scriptstyle\gamma$}
\enddc
\vspace*{0.3cm}

\begin{lemma}\label{branch locus}
The double cover $\pi: X\to \PSpc^3$ branches along the union of
the hyperplane arrangement $(H_1,H_2,\cdots,H_{8})$. \proof The
Galois group of $\pi$ is generated by $\iota_1$ in $G_{2}/G_{2}'$.
By the commutativity of the above diagram, the branch locus of
$\pi$ is the image of the fixed locus of $\iota_1$ under the
morphism $\gamma \circ q^3$. By Lemma \ref{in general position},
it is clear that the image is $\bigcup_{i=1}^{8}H_{i}$. \qed
\end{lemma}

By this lemma, the moduli of hyperelliptic curves of genus 3 are
embedded into the moduli space $\ModularCY$. We call the image
$\HyperCY$ the \emph{hyperelliptic locus}, which is
five-dimensional. In \cite{MSY}, the analogous submoduli were also
characterized as those six lines in general positions tangential
to a smooth conic of $\PSpc^2$, and it was shown that this
submoduli gives the family of Kummer surfaces. The Hodge structure
of CY threefold over the hyperelliptic locus is also special in
our case.

\begin{proposition}\label{wedge product of weight one HS}
Let $\tilde{X}$ be the canonical resolution of $X$. We have an
isomorphism of rational polarized Hodge structures
$$
\Cohom^3(\tilde{X},\Rationals) \cong \bigwedge^3
\Cohom^1(C,\Rationals).
$$
\proof As a consequence of Proposition \ref{XYRelation}, we know
that
$$
\Cohom^3(\tilde{X},\Rationals) \cong \Cohom^{3}(X,\Rationals).
$$
So it suffices to prove the isomorphism for $X$. For this purpose we
consider the following commutative diagram

\vspace*{0.5cm}\hspace*{3.5cm}
\begindc{\commdiag}[25]
\obj(1,3)[C3]{$C^3$} \obj(3,3)[S3C]{$S^3(C)$}
\obj(5,3)[JC]{$\Jacobian(C)$} \obj(3,1){$X$}
\mor{C3}{S3C}{$\scriptstyle\delta_1$}
\mor{S3C}{JC}{$\scriptstyle\varphi$}
\mor{C3}{$X$}{$\scriptstyle\delta$}
\mor{S3C}{$X$}{$\scriptstyle\delta_2$}
\enddc
\vspace*{0.3cm}

\noindent where $\varphi$ is the Abel-Jacobi map, $\delta_{1}$ is
the quotient map by the subgroup $S_3 \leq G_2'$ and $\delta_{2}$
is the projection map. One notes that, since $S_3$ is not normal
in $G_{2}'$, the map $\delta_{2}$ is only a finite morphism.
However, $\delta$ induces the embedding
$$
\delta^{*}: \Cohom^{3}(X,\Rationals) \cong
\Cohom^3(C^{3},\Rationals)^{G_2'} \hookrightarrow
\Cohom^{3}(C^{3},\Rationals).
$$
Since $\delta=\delta_{2}\circ \delta_{1}$, the pullback
$\delta_{2}^*$ gives the embedding
$$
H^3(X,\Rationals)\stackrel{\delta_{2}^*}{\longrightarrow}
H^3(S^3(C),\Rationals).
$$
By the Abel-Jacobi theorem, $\varphi$ is a birational morphism and
thus induces an isomorphism of Hodge structures on the middle
cohomology:
$$
\varphi^* : \Cohom^3(\Jacobian(C),\Rationals)
\stackrel{\cong}{\longrightarrow} \Cohom^3(S^3(C),\Rationals).
$$
In particular, $\dim_{\Rationals} \Cohom^3(S^3(C),\Rationals)=20$.
Since we computed before that the dimension of
$\Cohom^3(\tilde{X},\Rationals)$ is also $20$, the map
$\delta_{2}^*$ is in fact an isomorphism. The composition map
$$
\Cohom^{3}(X,\Rationals) \stackrel{\delta_{2}^*}{\longrightarrow}
\Cohom^3(S^3(C),\Rationals)
\stackrel{(\varphi^*)^{-1}}{\longrightarrow}
\Cohom^3(\Jacobian(C),\Rationals) \cong \bigwedge^{3}
\Cohom^{1}(C,\Rationals)
$$
gives the isomorphism required in the proposition.\qed
\end{proposition}

\begin{remark}
It is worthwhile to remark that the same construction and
arguments generalize to $n\geq 4$ cases. It will give a
$(2n-1)$-dimensional hyperelliptic locus in the
$n^{2}$-dimensional moduli of CY manifolds, over which the
primitive middle dimensional rational Hodge structures are wedge
products of weight 1 Hodge structures.
\end{remark}

\section{Characteristic Subvariety and Plethysm}
\label{sec:MethodDescr}

In this section, we will present two different methods to disprove
the modularity of $\ModularCY$. Our first method is to study a
series of invariants of IVHS, introduced in \cite{SZ}, which we call
\emph{characteristic subvarieties}. These invariants exploit the geometry
of the kernels of iterated Higgs fields of the associated system
of Hodge bundles with the given IVHS. In the case of Calabi-Yau
like PVHS over bounded symmetric domain, these invariants are
proved to be the \emph{characteristic bundles} introduced in
\cite{Mok} by N.\ Mok, which played a pivotal role in the proof of
the metric rigidity theorem of compact quotient of bounded symmetric
domains of rank $\geq 2$. The second method uses the idea of
\emph{plethysm} in representation theory (cf.\ \cite{FH}). For a
fixed simple complex Lie algebra $\mathfrak{g}$ the plethysm
describes the decompositions of representations derived
from a given irreducible representation of $\mathfrak{g}$.

\subsection{Characteristic Subvariety}

We first recall some results in \cite{SZ}. The bounded symmetric
domain
$$
D=\frac{\SpecialUnitary(3,3)}{\mathrm{S}(U(3)\times \Unitary(3))}
$$
is of rank 3. Let $\bb{W}$ be the Calabi-Yau like PVHS over
$\Gamma\backslash D$ and $(F,\eta)$ be the associated system of
Hodge bundles. By Theorem 3.3 in \cite{SZ} we have the following

\begin{lemma}\label{CharSubvariety}
For $k=1,2$ the $k$-th characteristic subvariety $\ca{S}_{k}$ of
$(F,\eta)$ coincides with $k$-th characteristic bundle. In
particular, for every point $x \in \Gamma\backslash D$,
$$
(\ca{S}_1)_{x} \cong \PSpc^2 \times \PSpc^2,
$$
and $(\ca{S}_2)_{x}$ is isomorphic to the determinantal
hypersurface in $\PSpc^8$.
\end{lemma}

Now we take a universal family $f: \ca{X} \rightarrow S$ of
$\ModularCY$. Let $\bb{V}:=R^{3}f_{*}\Rationals_\ca{X}$ and
$$
\left(
E=\bigoplus_{p+q=3}E^{p,q},\theta=\bigoplus_{p+q=3}\theta^{p,q}
\right)
$$
be the corresponding system of Hodge bundles. Since $\bb{V}$ is of
weight 3, we have also two characteristic subvarieties of
$(F,\eta)$, which are denoted by $\ca{R}_k$ for $k=1,2$. If $f$ is
a modular family, then the period mapping $\phi: S\hookrightarrow
\Gamma\backslash D$ will induce an isomorphism
$$
\phi^{*}\bb{W} \cong \bb{V},
$$
hence an isomorphism $\phi^{*}(E,\theta) \cong (F,\eta)$. This
implies the isomorphisms
$$
\phi^{*} \ca{S}_k \cong \ca{R}_k \qquad \text{for} \qquad k=1,2.
$$
Using Lemma \ref{CharSubvariety}, we then have the following

\begin{proposition}\label{CharSubProp}
If there exists a point $x\in S$ such that
$$
(\ca{R}_1)_{x} \not\cong \PSpc^2 \times \PSpc^2
$$
or $(\ca{R}_{2})_{x}$ is not isomorphic to the determinantal
hypersurface in $\PSpc^8$, then $f$ is not a modular family.
\end{proposition}

\begin{remark}
It was first pointed out by E.\ Viehweg that the iterated Higgs
fields for $(E,\theta)$ are surjective. Namely, the maps
$$
\theta^3 : S^{k} ( \ca{T}_{S} ) \longrightarrow
\Hom(E^{3,0},E^{3-k,k})
$$
are surjective for $1\leq k\leq 3$, where $\ca{T}_S$ denotes the
tangent bundle over $S$. If one of these maps were not surjective,
then the disproof of modularity of $\ModularCY$ would have been
obtained at this stage already. This phenomenon (or difficulty)
actually motivated the two latter authors to study the
characteristic subvariety in \cite{SZ}. It turned out that the
present work gives a non-trivial application of the theory of
characteristic subvarieties.
\end{remark}

\subsection{Plethysm}\label{PlethysmMethod}
The simple real Lie group $\SpecialUnitary(3,3)$ is a real form of
$\SpecialLinear(6,\Complex)$. By Weyl's unitary trick, one has an
equivalence of categories of finite dimensional complex
representations of $\SpecialUnitary(3,3)$ and finite dimensional
complex representations of
$\mathfrak{g}:=\SpecialLinearLie(6,\Complex)$. So the plethysm
problem for $\SpecialUnitary(3,3)$ is transformed into the
plethysm problem for $\mathfrak{g}$.\\

Let $V:=\Complex^6$ be the standard representation of
$\mathfrak{g}$. We shall study the plethysm for the fundamental
representation $W:=\bigwedge^3(V)$. In other words, we shall study
the decomposition of $\bb{S}_{\lambda}(W)$ for a Schur functor
$\bb{S}_{\lambda}$. The two simplest Schur functors are $S^{2}$
and $\bigwedge^{2}$. By Exercise 15.32 in \cite{FH} we have the
following decompositions:
\begin{equation}\label{S2W2Decomp}
S^2(W)=\Gamma_{10001}\oplus \Gamma_{00200},\quad\quad
\bigwedge^{2} W =\Gamma_{00000}\oplus \Gamma_{01010}
\end{equation}
By formula (15.17) in \cite{FH} it is easy to compute that
$$
\dim_\Complex \Gamma_{00000} = 1, \qquad \dim_\Complex
\Gamma_{10001} = 35, \qquad \dim_\Complex \Gamma_{01010} = 189
$$
and $\dim_\Complex \Gamma_{00200} = 175$. However, $\bigwedge^{2}
W$ will be of no use for us. That is because, considering $W$ as
$\SymplecticLie(20,\Complex)$-representation, one also has a
decomposition
$$
\bigwedge^{2} W = \Complex \oplus W'
$$
where $\Complex$ is the trivial representation of
$\SymplecticLie(20,\Complex)$ spanned by the symplectic form. On
the other hand, $S^2(W)$ is an irreducible representation of
$\SymplecticLie(20,\Complex)$. It is actually the adjoint
representation.

\begin{proposition}\label{PlethysmProp}
Let $f:\ca{X} \rightarrow S$ be a universal family of
$\ModularCY$. If $S^{2}(E,\theta)$ is not decomposed according to
the following pattern, then $f$ is not a modular family.
Explicitly,
$$
S^{2}(E,\theta)=(E_1,\theta_1)\oplus (E_2,\theta_2)
$$
where
\begin{eqnarray*}
  E_{1} &=&\phantom{E_{1}^{6,0}\oplus E_{1}^{5,1}\oplus}\ E_{1}^{4,2}\oplus
E_{1}^{3,3}\oplus
  E_{1}^{2,4}
  \phantom{E_{1}^{1,5}\oplus E_{1}^{0,6}}  \\
  E_{2} &=& E_{2}^{6,0}\oplus E_{2}^{5,1}\oplus E_{2}^{4,2} \oplus
E_{2}^{3,3} \oplus E_{2}^{2,4}\oplus E_{2}^{1,5}\oplus
E_{2}^{0,6}.
\end{eqnarray*}
Furthermore, the dimensions of Hodge bundles of $E_{1}$ are
respectively $0,0,9,17,9,0,0$ and those of $E_{2}$ are
$1,9,45,65,45,9,1$.

\proof The modularity of $f$ will imply a factorization of the
monodromy representation
$$
\rho: \pi_{1}(S)\longrightarrow
\SpecialUnitary(3,3)\stackrel{\bigwedge^3}{\longrightarrow}
\Symp(20,\Reals).
$$
Thus for any Schur functor $\bb{S}_{\lambda}$ the derived PVHS
$\bb{S}_{\lambda}(\bb{V})$ will decompose into irreducible
$\SpecialUnitary(3,3)$-representations. By the formula
\eqref{S2W2Decomp} and Deligne \cite{D1} Prop. 1.13, we have an
decomposition of PVHS
$$
S^{2}(\bb{V})=\bb{V}_1\oplus \bb{V}_2.
$$
The system of Hodge bundles $S^{2}(E,\theta)$ decomposes into
a direct sum of system of Hodge bundles accordingly,
$$
S^{2}(E,\theta)=(E_1,\theta_1)\oplus (E_2,\theta_2).
$$
Since $\bb{W}$ is of weight 3, $S^{2}(\bb{W})$ is of weight 6. One
can compute the Hodge numbers of $(E_i,\theta_i)$ for $i=1,2$ by
restricting the irreducible representations of
$\SpecialUnitary(3,3)$ to the center $\Unitary(1)$ of its maximal
compact subgroup $\mathrm{S}(\Unitary(3)\times \Unitary(3))$. If
$\IdMatrix$ denotes the $3 \times 3$-identity matrix, then
$$
\left\{ C_{z} := \begin{pmatrix} z \IdMatrix & 0 \\
0 & z^{-1} \IdMatrix \end{pmatrix} \in \GeneralLinear_6(\Complex)
\mid z\in \Unitary(1) \right\}
$$
is the center of $\mathrm{S}(\Unitary(3)\times \Unitary(3))$. We
choose the standard basis $(e_1,...,e_6)$ of $V = \Complex^6$ such
that
$$
C_z(e_i) = ze_i \quad \text{for} \quad 1 \leq i \leq 3 \qquad
\text{and} \qquad C_z(e_i) = z^{-1}e_i \quad \text{for} \quad 4
\leq i \leq 6.
$$
One notes that $\Gamma_{10001}$ is the unique nontrivial component
in $\Gamma_{10000}\otimes \Gamma_{00001}$. It is easy to compute
that $C_{z}$ acts on $\Gamma_{10001}$ with three characters
$z^2,z^0,z^{-2}$, and the dimensions of their eigenspaces are
respectively $9,17,9$. Then the characters of $C_{z}$ on the other
direct component $\Gamma_{00200}$ are
$z^6,z^4,z^{2},z^0,z^{-2},z^{-4},z^{-6}$, and their dimensions of
eigenspaces are computed to be $1,9,45,65,45,9,1$, respectively.
The proof of the proposition is complete.\qed
\end{proposition}

\section{The Jacobian Ring}
\label{sec:JacobiRing}

In the subsequent part we will carry out the strategies described
in section \ref{sec:MethodDescr} to the special family of CY
$3$-folds constructed in section \ref{sec:CYConstruction}. For this
purpose we let $S$ denote the moduli space of eight planes in $\PSpc^3$
in general positions. Every point $s \in S$ can be determined by a
matrix $A \in \Complex^{4 \times 8}$ with the property that all $(4
\times 4)$-minors of $A$ are non-zero. Furthermore, we let
$$
f : \tilde{\ca{X}} \longrightarrow S
$$
denote the universal family of $\ModularCY$ such that for every
every fiber $\tilde{X} := \tilde{\ca{X}}_s$ is obtained by
resolution of singularities from the ramified double cover $X
\rightarrow \PSpc^3$ associated to a certain matrix $A$ as described in
section \ref{sec:CYConstruction}. For our purposes it will be
necessary to give an explicit description of the PVHS $\bb{V} :=
R^3 f_* \Complex_\ca{X}$ and the
associated system $(E,\theta)$ of Higgs fields in every fiber.\\

First we give a description of $\bb{V}$ as a local system of
graded $\Complex$-vector spaces. Let $\sc{O}_S$ denote the sheaf
of holomorphic functions on $S$ and $a_{ij} \in
\Gamma(S,\sc{O}_S)$ the coordinate functions for $1 \leq i \leq 4$
and $1 \leq j \leq 8$. Furthermore, we let
$$
\sc{R} := \sc{O}_S[x_1,...,x_8,y_1,...,y_4]
$$
denote the free $\sc{O}_S$-algebra in $12$ indeterminates. For $p
\in \Naturals_0$ we define $\sc{R}^p$ to be the
$\sc{O}_S$-submodule of elements which have total degree $\deg_X =
2p$ in the variables $x_j$ and total degree $\deg_Y = p$ in the
variables $y_i$. We define a global sections $f_i,f \in
\Gamma(S,\sc{R})$ by $f_i := \sum_{j=1}^8 a_{ij} x_j^2$ for $1
\leq i \leq 4$ and $F := \sum_{i=1}^4 y_i f_i$. The twelve partial
derivatives
$$
\frac{\partial F}{\partial x_j} \quad \text{for} \quad 1 \leq j
\leq 8 \qquad \text{and} \qquad \frac{\partial F}{\partial y_i}
\quad \text{for} \quad 1 \leq i \leq 4
$$
generate an ideal sheaf in $\sc{R}$ which we denote by $\sc{I}$.
Finally, we let the group $G_1$ from section
\ref{sec:CYConstruction} act on the sheaf $\sc{R}$ by sending $a =
(a_1,...,a_8) \mapsto \sigma_a$ with
$$
\sigma_a(x_i) = (-1)^{a_i} x_i \quad \text{for $1 \leq i \leq 8$}
\qquad \text{and} \qquad \sigma_a(y_j) = y_j \quad \text{for $1
\leq j \leq 4$}.
$$
Then obviously $\sigma_a(\sc{I}) \subseteq \sc{I}$ holds for all
$a \in G_1$. Now we obtain the following explicit description of
our PVHS $\bb{V}$.

\begin{proposition}\label{JacRingCSpc}
There is a canonical isomorphism of local systems
$$
\bb{V} \otimes_\Complex \sc{O}_S \cong (\sc{R}/\sc{I})^{N_1}
$$
which maps $\bb{V}^{3-p,p} \otimes_\Complex \sc{O}_S$ onto the
submodule generated by $\sc{R}^p$ for $0 \leq p \leq 3$.

\proof Let $g : \ca{Y} \rightarrow S$ denote the family of
intersections of four quadrics in $\PSpc^7$ as constructed in
section \ref{sec:CYConstruction}, i.e.\ for every $s = A =
(a_{ij})\in S$ the fiber $\ca{Y}_s$ in the intersection of
quadrics given by the equations \eqref{FourQuadrics}. Furthermore,
by $\bb{W}:=R^3 g_* \Complex_\ca{Y}$ we denote the associated
PVHS. According to Proposition \ref{XYRelation} we have a
canonical isomorphism $\bb{V} \cong \bb{W}^{N_1}$, so that it
remains to establish the isomorphism
\begin{equation}\label{FixedIso}
\bb{W}^{N_1} \otimes_\Complex \sc{O}_S \cong
(\sc{R}/\sc{I})^{N_1}.
\end{equation}
First we show that $\bb{W} \otimes_\Complex \sc{O}_S \cong
\sc{R}/\sc{I}$. This is a special case of Proposition 2.2.10 in
\cite{Nagel}, and although it is stated only for individual
varieties, the result carries over to algebraic families. Here we
just sketch the essential steps. Let $\PSpc^7_S$ denote the
projective $7$-space over $S$ on which the coherent sheaf
$$
\sc{E} := \sc{O}_{\PSpc^7_S}(2)^{\oplus 4}
$$
is defined, and let $P := \bb{P}(\sc{E})$ denote the associated
projective bundle. Then $P$ contains a toric hypersurface
$\hat{\ca{Y}}$ given by the equation $F = \sum_{i=1}^4 y_i f_i$
from above. Let $\pi : P \rightarrow \PSpc^7_S$ denote the
canonical projection, extend $g$ to a map $g : \PSpc^7_S
\rightarrow S$ and let $h := g \circ \pi$. Then the embedding
$\pi^{-1}(\ca{Y}) \hookrightarrow \hat{\ca{Y}}$ induces a natural
isomorphism
$$
R^9 h_* \Complex_{\hat{\ca{Y}}} \cong \bb{W} \otimes
\Cohom^6(\PSpc^3,\Complex)
$$
of PVHS on $S$, the right part of the tensor product being
constant of rank
one.\\

Now let $\ca{V} := P \setminus \ca{Y}$ denote the open complement
of $\ca{Y}$. Then the Gysin sequence relating the PVHS's of $P$,
$\hat{\ca{Y}}$ and $\ca{V}$ gives rise to an isomorphism $R^9 h_*
\Complex_{\hat{\ca{Y}}} \cong R^{10} h_* \Complex_{\ca{V}}$ of
PVHS. In order to compute the latter, we make use of de Rham's
theorem which enables us to describe the cohomology
$$
R^{10} h_* \Complex_{\ca{V}} \otimes_\Complex \sc{O}_S \cong
\bbm{R}^{10} h_* \Omega^\cdot_{\ca{V}|S} \cong \bbm{R}^{10} h_*
\Omega^\cdot_{P|S}(*\hat{\ca{Y}})
$$
in terms of the sheaf $\Omega^\cdot_{P|S}(*\hat{\ca{Y}})$ of
relative differentials on $P$ with poles along $\hat{\ca{Y}}$,
where the functor $\bbm{R}^{10} h_*$ denotes hypercohomology.
Since the sheaves $\Omega^i_{P|S}(m \hat{\ca{Y}})$ are acyclic for
$i,m > 0$, it can be computed by taking global sections. That is,
if $\sc{Z}_{P|S} \subseteq \Omega^{10}_{P|S}(*\hat{\ca{Y}})$
denotes the subsheaf of closed differentials and $\sc{B}_{P|S}$
the subsheaf of exact ones, then simply
$$
\bbm{R}^{10} h_* \Omega^\cdot_{P|S}(*\hat{\ca{Y}}) \cong (h_*
\sc{Z}_{P|S}) / (h_* \sc{B}_{P|S}).
$$
The sections of $h_* \Omega^{10}_{P|S}(*\hat{\ca{Y}})$ can be
described in terms of the $\sc{O}_S$-algebra $\sc{R}$. Namely, let
$\omega_0$ denote the homogeneous differential form
$$
\hat{\Omega} := dx_1 \wedge \cdots \wedge dx_8 \wedge dy_1 \wedge
\cdots \wedge dy_4
$$
and define the vector fields $\vartheta_i:=\partial/\partial x_i$
and $\lambda_j:=\partial/\partial y_j$ for $1 \leq i \leq 8$ and
$1 \leq j \leq 4$. If we put
$$
\theta_1 := \sum_{j=1}^4 y_j \lambda_j \qquad , \qquad \theta_2 :=
\sum_{i=1}^8 x_i \vartheta_i - 2\sum_{j=1}^4 y_j \lambda_j
$$
and $\Omega := \theta_1\theta_2(\hat{\Omega})$, then every section
$\omega$ of $h_* \Omega^\cdot_{P|S}(p\hat{\ca{Y}})$ can be written
in the form
\begin{equation}\label{GenDiff}
\omega = \frac{H\Omega}{F^{p+4}} \qquad \text{where $H$ is a
section of $\sc{R}^p$}.
\end{equation}
In degree $9$, any section $\psi$ of $h_*
\Omega^{9}_{P|S}(*\hat{\ca{Y}})$ can be written as
$$
\psi = \frac{\sum_{i=1}^8 G_i \Omega_i - \sum_{j=1}^4 H_j
\Omega_j'}{F^{p+4}}
$$
where $\Omega_i := \theta_1\theta_2\vartheta_i$, $\Omega_j' :=
\theta_1\theta_2\lambda_j$ and $G_i,H_j$ are sections of $\sc{R}$
such that $\deg_X(G_i) = 2p+1,\deg_Y(G_i) = p$ and $\deg_X(H_j) =
2p$, $\deg_Y(H_j) = p+1$ for $1 \leq i \leq 8$ and $1 \leq j \leq
4$. Its exterior derivative is $d\psi = H \Omega / F^{p+5}$ where
$$
H = 2 \sum_{i=1}^8 \frac{\partial F}{\partial x_i} G_i + 2
\sum_{j=1}^4 \frac{\partial F}{\partial y_j} h_j - F \left(
\sum_{i=1}^8 \frac{\partial G_i}{\partial x_i} + \sum_{j=1}^4
\frac{\partial H_j}{\partial y_j} \right).
$$
We see that $\omega$ can be reduced to lower pole order if and
only if the section $h$ is a section of the ideal sheaf $\sc{I}$.
This shows that
$$
(h_* \sc{Z}_{P|S}) / (h_* \sc{B}_{P|S}) \cong \sc{R}/\sc{I}.
$$
Combing all isomorphisms, the desired assertion $\bb{W}
\otimes_\Complex \sc{O}_S \cong \sc{R}/\sc{I}$ follows. Observing
that the action of $G_1$ on $\bb{W} \otimes_\Complex \sc{O}_S$ is
compatible with the action defined above on $\sc{R}/\sc{I}$, we
obtain $\bb{W}^{N_1} \otimes_\Complex
\sc{O}_S \cong (\sc{R}/\sc{I})^{N_1}$.\\

In order to prove the refined statement on the grading, notice
that by the above construction
$$
\bb{W}^{3-p,p} \cong R^{6-p,3+p} h_* \Complex_{\hat{\ca{Y}}} \cong
R^{7-p,3+p} h_* \Complex_{\ca{V}}.
$$
By the comparison of Hodge and pole filtration, the part $(F^{3+p}
h_* \Complex_{\ca{V}}) \otimes_\Complex \sc{O}_S$ coincides with
the subsheaf of $\bbm{R}^{10} h_*
\Omega^\cdot_{P|S}(*\hat{\ca{Y}})$ generated by differentials of
pole order $\geq p+4$. This shows that $\bb{W}^{3-p,p}
\otimes_\Complex \sc{O}_S$ corresponds to the subsheaf of
$\sc{R}/\sc{I}$ generated by $\sc{R}^p$.\qed
\end{proposition}

The description of the local system $\bb{V}$ in terms of the
Jacobian ring $\sc{R}/\sc{I}$ admits an explicit computation the
Gauss-Manin connection and the Higgs field in one-parameter
families. Let $h : \hat{\ca{Y}} \rightarrow S$ denote the family
of toric hypersurfaces that we used in the proof of Proposition
\ref{JacRingCSpc}. Furthermore, let $T$ denote an open subset of
$\ASpc^1$ and $h : \hat{\ca{Y}}_T \rightarrow T$ be the family
obtained by restriction. Over $T$ the defining equation of
$\hat{\ca{Y}}$ inside the toric variety $P$ is given by an
equation $F = 0$ with $F \in \Complex(t)[x_1,...,y_4]$, and the
Gauss-Manin connection
$$
\nabla : R^9 h_* \Complex_{\hat{\ca{Y}}_T} \longrightarrow R^9 h_*
\Complex_{\hat{\ca{Y}}_T} \otimes_{\sc{O}_T} \Omega^1_T
$$
acts on the de Rham cohomology of the complement by
\begin{equation}\label{GMMap}
\omega = \frac{H\Omega}{F^{p+4}} \mapsto -(p+4) (\partial_t F)
\frac{H\Omega}{F^{p+5}} \otimes dt
\end{equation}
provided that the section $H$ of $\sc{R}_T$ is chosen such that
$\partial_t H = 0$. Thus if one fixes a local basis of
$(\sc{R}_T/\sc{I}_T)^{N_1}$ given by polynomials over the function
field $\Complex(t)$, one can compute a representation matrix of
$\nabla$ by applying the map \eqref{GMMap} to all basis elements
and reducing them with respect to the basis. A representation
matrix for the Higgs field
$$
\theta : R^9 h_* \Complex_{\hat{\ca{Y}}_T} \longrightarrow R^9 h_*
\Complex_{\hat{\ca{Y}}_T} \otimes_{\sc{O}_T} \Omega^1_T
$$
is obtained by projecting the images of the basis elements inside
the $\sc{R}_T^p$-part onto the subspace generated by
$\sc{R}_T^{p+1}$, for $0 \leq p \leq 3$. By the canonical
isomorphism of Proposition \ref{JacRingCSpc} this also yields a
local representation matrix of $\theta : \bb{V}_T \rightarrow
\bb{V}_T \otimes_{\sc{O}_T} \Omega^1_T$, or equivalently, of
\begin{equation}\label{LocalHiggsField}
\theta : \ca{T}_T \otimes_{\sc{O}_T} \bb{V}_T \longrightarrow
\bb{V}_T.
\end{equation}
For our purposes it will be sufficient to compute the map
\eqref{LocalHiggsField} in the infinitesimal neighborhood of a
point $x \in T$, which turns out to be much easier. Let $\tilde{X}
:= \tilde{\ca{X}}_x$ denote the fiber at $x$. We have an exact
sequence of vector bundles over $\tilde{X}$ given by
\begin{equation}\label{LongExactTangent}
0 \longrightarrow \ca{T}_{\tilde{X}} \longrightarrow
\ca{T}_{\tilde{\ca{X}}|\tilde{X}} \longrightarrow f^*( \ca{T}_S
)_{|\tilde{X}} \longrightarrow 0
\end{equation}
where the vertial bars mean restriction. The bundle on the right
hand side is trival with generic fiber $T_{S,x}$, the tangent
space of $S$ at $x$. Since $\tilde{X}$ is compact, all sections of
the trivial bundle are constant, so that
$$
T_{S,x} = \Cohom^0(X,f^*( \ca{T}_S )_{|\tilde{X}})
$$
holds. Now the long exact cohomology sequence associated to the
short exact sequence (\ref{LongExactTangent}) yields a map
$$
\rho : T_{S,x} \longrightarrow
\Cohom^1(\tilde{X},\ca{T}_{\tilde{X}}) \quad,
$$
the {\em Kodaira-Spencer map}. It is know to be an isomorphism.\\

Let $R$ denote the stalk of the local system
$(\sc{R}/\sc{I})^{N_1}$ at $x$, and by $R_p$ the stalks of the
images of $\sc{R}_p$, for $0 \leq p \leq 3$. Then $R =
\bigoplus_{p=0}^3 R_p$ is a finite-dimensional $\Complex$-algebra.

\begin{lemma} \label{KodairaLemma}
There is a canonical isomorphism $R_1 \cong
\Cohom^1(\tilde{X},\ca{T}_{\tilde{X}})$. \proof Since $\tilde{X}$
is a Calabi-Yau manifold, the canonical bundle $\sc{K}_{\tilde{X}}
= \Omega^3_{\tilde{X}}$ is trivial, which gives rise to a natural
identification
$$
\ca{T}_{\tilde{X}} = (\Omega^1_{\tilde{X}})^* \cong
\Omega^2_{\tilde{X}}.
$$
It implies that $\Cohom^1(\tilde{X},\ca{T}_{\tilde{X}})$ is
isomorphic to $\Cohom^{2,1}(\tilde{X}) =
\Cohom^1(\tilde{X},\Omega^2_{\tilde{X}})$. On the other hand, if
we specialize the isomorphism from Proposition \ref{JacRingCSpc}
to the stalks at $x$, we obtain $\Cohom^{2,1}(\tilde{X}) \cong
R_1$.\qed
\end{lemma}

\begin{proposition} \label{ThetaMapMultMap}
For $0 \leq p \leq 2$ there is a commutative diagram

\vspace*{0.5cm}\hspace*{2.0cm}
\begindc{\commdiag}[25]
\obj(1,3)[TH]{\hspace*{1.0cm}$T_{S,x} \otimes
\Cohom^{3-p,p}(\tilde{X})$}
\obj(6,3)[H]{$\Cohom^{2-p,p+1}(\tilde{X})$} \obj(1,1)[RR]{$R_1
\otimes R_p$} \obj(6,1)[R]{$R_{p+1}$}
\mor{TH}{H}{$\scriptstyle\theta$}
\mor{TH}{RR}{$\scriptstyle\cong$} \mor{H}{R}{$\scriptstyle\cong$}
\mor{RR}{R}{$\scriptstyle\mu$}
\enddc
\vspace*{0.3cm}

\noindent where the vertical arrows are induced by the
Kodaira-Spencer map and Proposition \ref{JacRingCSpc}, and where
the lower horizontal arrow denotes multiplication on the graded
$\Complex$-algebra $R$. \proof It is known that the derivation of
a cohomology class in $\Cohom^{3-p,p}(\tilde{X})$ with respect to
a tangent direction $v \in T_{S,x}$ is given by the cup product
$$
\Cohom^1(\tilde{X},\ca{T}_{\tilde{X}}) \otimes
\Cohom^q(\tilde{X},\Omega^p_{\tilde{X}})
\stackrel{\cup}{\longrightarrow}
\Cohom^{q+1}(\tilde{X},\Omega^{p-1}_{\tilde{X}})
$$
with the Kodaira-Spencer class $\rho(v)$ (see e.g.\
\cite{Periods}, Lemma 5.3.3). In the de Rham cohomology of the
toric hypersurface $\hat{\ca{Y}}$, the cup product between
cohomology classes corresponds to the wedge product between
differential forms. Furthermore, we have seen in \eqref{GenDiff}
that every differential is defined by a polynomial in $R$. It can
be checked easily that the multiplication of polynomials
corresponds to the wedge product of the corresponding differential
forms.\qed
\end{proposition}

For later use we need an explicit, fiberwise description of the
characteristic subvarieties $\ca{R}_k$ introduced in section
\ref{sec:MethodDescr} associated to our special universal family
$f : \tilde{\ca{X}} \rightarrow S$. To this end we introduce the
symmetric algebra $S^\cdot(R_1^*)$ over the dual of $R_1$, which
is the homogeneous coordinate ring of $\PSpc(R_1^*)$. Taking the
multiplication map to its dual, we obtain a linear map
$$
\mu^* : R_2^* \longrightarrow S^2(R_1^*) \quad,
$$
and we let $\fr{a}_1$ denote the ideal generated by the image of
$\mu^*$. Similiarly, we let $\fr{a}_2$ denote the ideal generated
by the image of the dualized multiplication map $\mu^* : R_3^*
\rightarrow S^3(R_1^*)$. Then the fibers of the characteristc
varieties are obtained in the following way.

\begin{lemma} \label{CharSubvComp}
For a point $x \in S$ as above and $k=1,2$, the fiber of the
$k$-th characteristic subvariety $(\ca{R}_k)_x$ is isomorphic to
the projective subvariety $Z_k := \Proj(A_k)$ of $\PSpc(R_1^*)$,
where $A_k$ denotes the graded quotient ring
$S^\cdot(R_1^*)/\fr{a}_k$. \proof We recall the definition of the
$k$-th characteristic subvariety as given in \cite{SZ}. For our
system $(E,\theta)$ of Hodge bundles, the $(k+1)$-st iterated
Higgs field defines a map
$$
\theta^{k+1} : S^{k+1}(\ca{T}_S) \longrightarrow
\Hom(E^{3,0},E^{2-k,k+1})
$$
whose kernel we denote by $\sc{I}_k$. Then $\ca{R}_k =
\Proj(\sc{I}_k)$ as a subvariety of $\PSpc(\ca{T}_S)$. For $k=1$
the stalk $(\sc{I}_1)_x$ at $x \in S$ is the kernel of
$$
\theta^2_x : S^2(T_{S,x}) \cong S^2(R_1) \longrightarrow
\Hom(R_3,R_2) \cong R_2 \quad,
$$
the first isomorphism coming from Lemma \ref{KodairaLemma} and the
Kodaira-Spencer map, the second being a consequence of the fact
that $R_3$ is one-dimensional. If we dualize this map, up to a
non-zero constant we obtain $\mu^*$, and the kernel of
$\theta^2_x$ is isomorphic to $S^2(R_1^*)/\fr{a}_1$, the cokernel
of $\mu^*$. Since this quotient generates $A_1$, we obtain $Z_1
\cong (\ca{R}_1)_x$. The proof for $k=2$ is similar.\qed
\end{lemma}

\section{Proofs of the Main Theorems}
\label{sec:Proofs}

We recall some basic notions from computational commutative
algebra. Let $K$ be a field and $R:=K[x_1,...,x_n]$ the polynomial
ring in $n$ indeterminates. A {\em monomial ordering} is a total
ordering $\prec$ on the set of monomials in $R$ such that $f \prec
g$ implies $fh \prec gh$ for monomials $f,g,h \in R$. In our
computations we will use the {\em graded lexicographical
ordering}, which is defined as follows: First one fixes an
ordering on the set of indeterminates by requiring $x_1 \succ x_2
\succ \cdots \succ x_n$. Now let
$$
f = y_1 y_2 \cdots y_r \qquad \text{and} \qquad g = z_1 z_2 \cdots
z_s \qquad \text{with} \quad y_i,z_i \in \{ x_1,...,x_n \} \quad
\text{for all $i$}
$$
such that $y_i \succ y_j$ or $y_i = y_j$ for $i \leq j$, and
similary for the factors of $g$. Then by definition $f \succ g$ if
either $r > s$ or $r = s$ and there is an $m \in \Naturals_0$ such
that $y_i = z_i$ for $1 \leq i \leq m$ and
$y_{m+1} \succ z_{m+1}$.\\

The total ordering on the monomials extends to a partial ordering
on $R$ by defining $f \prec g$ iff the maximal monomial of $f$ is
smaller than the maximal monomial of $g$. Furthermore, zero is
defined to be the least element in $R$. If $\fr{a} \subseteq R$ is
an ideal, then we say that an element $f \in R$ is in {\em normal
form} with respect to $\fr{a}$ and write $f = \NormalForm(f)$ if
$f$ is minimal inside the coset $f + \fr{a}$. It can be shown that
the normal form is unique; in particular, $\NormalForm(f) = 0$ if
and only
if $f \in \fr{a}$.\\

Let $f : \tilde{\ca{X}} \longrightarrow S$ denote the family of
CYs defined at the beginning of section \ref{sec:JacobiRing}. In
order to prove the theorems from section \ref{sec:Intro}, it
suffices to consider one particular fiber of this family. Let
$\lambda_j := j$ for $1 \leq j \leq 8$ and define the matrix $A
\in \Complex^{4,8}$ by $a_{ij} := \lambda_j^i$ for $1 \leq i \leq
4$ and $1 \leq j \leq 8$. We define $x_0 \in S$ to be the point
corresponding to the matrix $A$ and let $\tilde{X} :=
\tilde{\ca{X}}_{x_0}$ denote its fiber. For $p=0,...,3$ let $R_p$
be the ring defined before Lemma \ref{KodairaLemma}. By
Proposition \ref{JacRingCSpc}, $R_p^{N_1}$ is isomorphic to
$H^{3-p,p}(\tilde{X})$ for $0 \leq p \leq 3$.

\begin{lemma}\label{BasisLemma}
The following elements constitute a basis of $R_p^{N_1}$ for
$p=0,...,3$.
$$
\begin{array}{ll}
p = 0 & 1 \\
p = 1 & x_5^2 y_2,~ x_5^2y_3,~ x_5^2y_4,~ x_6^2y_2,~ x_6^2y_3,~
x_6^2y_4,~
x_7^2y_2,~ x_7^2y_3,~ x_7^2y_4 \\
p = 2 & x_6^4y_3^2,~ x_6^4y_3y_4,~ x_6^4y_4^2,~ x_6^2x_7^2y_3^2,~
x_6^2x_7^2y_3y_4,~ x_6^2x_7^2y_4^2,~ x_7^4y_3^2,~ x_7^4y_3y_4,~ x_7^4y_4^2 \\
p = 3 & x_7^6y_4^3
\end{array}
$$
\proof For each $p$ we list all monomials with $\deg_X = 2p$ and
$\deg_Y = p$. If we let $e_1,...,e_8$ denote the canonical basis
of $(\Integers/2\Integers)^8$, then $N_1$ is generated by the set
$$
B := \{ e_i + e_{i+1} ~\mid~ 1 \leq i \leq 7 \} \cup \{ e_1 + e_8
\}.
$$
We remove all elements $g$ from the list with $\sigma_a(g) \neq g$
for some $a \in B$ or with $\NormalForm(g) \neq g$. By uniqueness
and linearity of the normal form, the remaining elements are
linearly independent in $R_p^{N_1}$. Since the Betti numbers of
$\tilde{X}$ are $1,9,9,1$, respectively, the assertion
follows.\qed
\end{lemma}

\begin{proposition} \label{CharSubvThm}
The fiber $(\ca{R}_1)_{x_0}$ of the first characteristic
subvariety at $x_0$ is two-dimensional. \proof By Lemma
\ref{CharSubvComp} we have to compute the ideal $\fr{a}_1 \subset
S^\cdot(R_1^*)$ which is generated by the image of the dual
multiplication map $\mu^* : R_2^* \rightarrow S^2(R_1^*)$. Let
$v_1,...,v_9$ denote the basis of $R_1$ and $w_1,...,w_9$ the
basis of $R_2$ as defined in \ref{BasisLemma}. Furthermore, we fix
a bijection
$$
\varphi : \{ (i,j) \in \Naturals^2 ~\mid~ 1 \leq i \leq j \leq 9
\} \stackrel{\sim}{\longrightarrow} \{ 1,...,45 \}
$$
and put $u_{\varphi(i,j)} := v_i v_j$. The first step is to
compute a representation matrix of the multiplication map
$$
\mu : S^2(R_1) \longrightarrow R_2
$$
with respect to the basis $u_1,...,u_{45}$ and $w_1,...,w_9$. By
computing the normal forms of elements with respect to the
Jacobian ideal $\sc{I}_{x_0} \subseteq R$, we determine
$c_{\varphi(i,j)k} \in \Rationals$ such that
$$
\NormalForm(v_i v_j) = \sum_{k=1}^9 c_{\varphi(i,j)k} w_k \qquad
\text{for} \quad 1 \leq i,j \leq 9.
$$
Then $C := (c_{\ell k}) \in \Complex^{45,9}$ is the desired
representation matrix. Its transpose represents $\mu^*$ with
respect to the dual basis $w_1^*,...,w_9^*$ and
$u_1^*,...,u_{45}^*$.\medskip

Notice that $(v_i v_j)^* = 2 v_i^* v_j^*$ for $1 \leq i,j \leq 9$.
Thus if we define
$$
\tilde{c}_{k\ell} :=
\begin{cases} c_{k\ell} & k = \varphi(i,j), \quad i = j \\
2c_{k\ell} & k = \varphi(i,j), \quad i \neq j \end{cases}
$$
then $\transpose{\tilde{C}} \in \Complex^{9,45}$ is a
representation matrix of $\mu^*$ with respect to $w_1^*,...,w_9^*$
and $\tilde{u}_1,...,\tilde{u}_{45}$, where
$\tilde{u}_{\varphi(i,j)} := v_i^* v_j^*$. Each row corresponds to
one generator of $\fr{a}_1$ in $S^2(R_1^*)$. Furthermore, the
choice of a basis $v_1^*,...,v_9^*$ admits a natural identification
$$
\Proj(S^\cdot(R_1^*)) \cong \PSpc^8.
$$
Let $z_1,...,z_9$ denote a new set of indeterminates. If we define
$f_1,...,f_9$ by
$$
f_\ell := \sum_{i=1}^9 \sum_{j=i}^9 \tilde{c}_{\varphi(i,j)\ell}
z_i z_j
$$
then the variety in $\PSpc^8$ defined by $f_1 = \cdots = f_9 = 0$
is isomorphic to $(\ca{R}_1)_{x_0}$. Since the matrix $\tilde{C}$
is known in explicit term, we can use computer algebra to compute
its dimension. We obtain $\dim (\ca{R}_1)_{x_0} = 2$.\qed
\end{proposition}

In order to prove the non-modularity of $f$, by Proposition
\ref{CharSubProp} it is sufficient to determine a single points $x
\in S$ such that the fiber $(\ca{R}_1)_x$ is not isomorphic to
$\PSpc^2 \times \PSpc^2$. By Proposition \ref{CharSubvComp}, the
fiber $(\ca{R}_1)_{x_0}$ is only two-dimensional. Thus both $f$
and $\ModularCY$ cannot be modular, and Theorem
\ref{MainTheorems} is proved.\\

\textbf{Proof of Theorem \ref{NonExtHyperLocus}:} The proof will
be achieved by contradiction. Let $\HyperCY\subset \HyperCY'\subset \ModularCY$
be an extension as described in the theorem, and let $f:\mathfrak{X}\to S$
denote a universal
family. Let $Z$ be the sublocus in $S$ mapping to $\HyperCY'$, and
$g=f|_{Z}: \mathfrak{X}|_{f^{-1}(Z)}\to Z$ be the corresponding
subfamily. Then $g$ is a modular family with respect to
$(\Symp(6,\Reals)/\Unitary(3),\bb W)$ where $\bb{W}$ is the
Calabi-Yau like PVHS over $\Symp(6,\Reals)/\Unitary(3)$(cf.
\cite{SZ}). Let $(F,\eta)$ be the corresponding Higgs bundle of
the subfamily $g$. By Theorem 3.3 in \cite{SZ}, the Higgs bundle
$(F,\eta)$ has two characteristic subvarieties and the fibers of
the first characteristic subvariety are all isomorphic to
$\PSpc^2$. Take one point $x\in Z$, and denote by
$(\ca{R}'_{1})_x$ be the fiber over $x$ of the first
characteristic subvariety of $(F,\eta)$. By the geometric
description of the characteristic subvariety in Lemma 3.2 \cite{SZ},
we know that in $\bb{P}(T_{S,x})$ the equality
$$
(\ca{R}'_{1})_x=(\ca{R}_{1})_x\cap \bb{P}(T_{Z,x})
$$
holds. Now if $\dim(\ca{R}_{1})_x=2$, then we neccessarily have an isomorphism
$$
(\ca{R}_{1})_x=(\ca{R}'_{1})_x\simeq \PSpc^2.
$$
Since our computation is local, we simply take the point $x$ to be
the same point as used in the above proof of Theorem
\ref{MainTheorems}. The arithmetic genus of $(\ca{R}_{1})_x$ is
calculated to be -41, whereas the arithmetic genus of $\PSpc^2$ is
0. So $(\ca{R}_{1})_x$ is non-isomorphic to $\PSpc^2$. Therefore
such an extension does not exist. \qed
\\

Now we give a second proof of Theorem \ref{MainTheorems} which is
based on the plethysm method described in subsection
\ref{PlethysmMethod}. As before by $(E,\theta)$ we denote the
Hodge bundle associated to the family $f : \tilde{X} \rightarrow
S$. The Higgs field
$$
\theta_{x_0} : T_{S,x_0} \otimes E^{3,0}_{x_0} \longrightarrow
E^{2,1}_{x_0}
$$
induces in a natural way a linear map
$$
S^2(\theta_{x_0}) : T_{S,x_0} \otimes S^2(E_{x_0})^{6,0}
\longrightarrow S^2(E_{x_0})^{5,1}
$$
on the symmetric $2$-space. By threefold iteration we obtain
$$
S^2(\theta_{x_0}^3) : S^3(T_{S,x_0}) \otimes S^2(E_{x_0})^{6,0}
\longrightarrow S^2(E_{x_0})^{3,3}.
$$

\begin{proposition}\label{PlethysmCompThm}
The image of $S^2(\theta_{x_0}^3)$ is $78$-dimensional. \proof By
Proposition \ref{ThetaMapMultMap} it is sufficient to compute the
image of the linear map
$$
S^2(\mu^3) : S^3(R_1) \otimes S^2(R_0) \longrightarrow S^2(R)
$$
induced by the multiplication map $\mu : R_1 \otimes R_0
\rightarrow R_1$. Let $v_1,...,v_{20}$ denote the basis of $R =
\oplus_{p=0}^3 R_p$ specified in Lemma \ref{BasisLemma}; in
particular, we assume that the elements $w_k := v_{k+1}$ span the subspace
$R_1$, where $1 \leq k \leq 9$. For each $k$ we determine the
coefficients $c_{ij}^{(k)} \in \Complex$ such that
$$
v_i w_k = \sum_{j=1}^{20} c_{ij}^{(k)} v_j \qquad \text{for $1
\leq i \leq 20$ and $1 \leq j \leq 9$.}
$$
This is done by reduction modulo $\sc{I}_{x_0}$ as in the proof of
Theorem \ref{CharSubvThm}. Then $C_k := c_{ij}^{(k)}$ is a
representation matrix of the linear map
$$
\mu_{w_k} : R \longrightarrow R, \qquad v \mapsto v w_k
$$
with respect to $v_1,...,v_{20}$. With these matrices it now an
easy task to compute the induced action of $\mu_{w_k}$ on
$S^2(R)$. Fix a bijection
$$
\varphi : \{ (i,j) \in \Naturals^2 ~\mid~ 1 \leq i \leq j \leq 20
\} \stackrel{\sim}{\longrightarrow} \{ 1,...,210 \}
$$
and define a basis $u_1,...,u_{210}$ of $S^2(R)$ by
$u_{\varphi(i,j)} := v_i v_j$ for $1 \leq i \leq j \leq 20$. Then
$S^2(\mu_{w_k})$ acts on this basis by
\begin{equation}\label{Symm2Formula}
S^2(\mu_{w_k})(u_{\varphi(i,j)}) := \sum_{\ell=1}^{20}
c_{i\ell}^{(k)} v_j v_{\ell} + \sum_{\ell=1}^{20} c_{j\ell}^{(k)}
v_i v_{\ell}.
\end{equation}
The subspace $U^{6,0}$ of $S^2(R)$ of degree zero is
one-dimensional and generated by $u_{\varphi(0,0)}$. Now the space
$S^2(\mu)(U^{6,0})$ is generated by the images of all maps
$S^2(\mu_{w_k})$ ($k=1,...,9$) applied to $u_{\varphi(0,0)}$. By
\eqref{Symm2Formula} and computational linear algebra it turns out
to be $9$-dimensional and thus all of $U^{5,1}$. Applying all maps
$S^2(\mu_{w_k})$ to $U^{5,1}$ we obtain a subspace of $S^2(R)$ of
dimension $45$ contained in $U^{4,2}$, and a third application
yields a $78$-dimensional subspace of $U^{3,3}$.\qed
\end{proposition}

Now we explain how Proposition \ref{PlethysmCompThm} implies
Theorem \ref{MainTheorems}. If $f : \tilde{\ca{X}} \rightarrow S$
were modular, by Proposition \ref{PlethysmProp} there would be a
decomposition of Hodge bundles
$$
S^2(E,\theta) = (E_1,\theta_1) \oplus (E_2,\theta_2) \quad,
$$
such that each graded piece $E_i^{6-p,p}$ has a specific
dimension. This decomposition would exist in any fiber. In
particular, the image of $S^2(E_{2,x_0})^{6,0} =
S^2(E_{x_0})^{6,0}$ under the iterated Higgs field
$$
S^2(\theta_{x_0}^3) : S^3(T_{S,x_0}) \otimes S^2(E_{x_0})^{6,0}
\longrightarrow S^2(E_{x_0})^{3,3}
$$
would be contained in $S^2(E_{2,x_0})^{3,3}$ and thus be at most
$65$-dimensional. But since the image of $S^2(\theta_{x_0}^3)$ has
dimension $78$, the decomposition cannot exist.\\

\textbf{Proof of Theorem \ref{ZarClosMon}:} Let
$\rho: \pi_{1}(S)\to \SymplecticGroup(20,\Reals)$
be the monodromy representation. We know that $\rho$ is
irreducible since the VHS $\bb{V}$ is irreducible. Let $G$ be the
Zariski closure of the monodromy group in $\SymplecticGroup(20,\Reals)$.
Then we have a factorization:
$$
\rho: \pi_{1}(S)\longrightarrow
G\stackrel{\varrho}{\longrightarrow} \SymplecticGroup(20,\Reals),
$$
and $\varrho: G\to \SymplecticGroup(20,\Reals)$ is irreducible. If $G$ is not
the whole group, $G$ must be a proper Lie subgroup of
$\SymplecticGroup(20,\Reals)$. We will now derive a contradiction by a sequence
of steps.\\

\noindent
{\em Step 1.} Differentiating $\varrho$ we pass to the real Lie algebra
monomorphism
$$
\chi: \mathfrak{g} \longrightarrow \SymplecticLie(20,\Reals)
$$
where $\mathfrak{g}=\Lie(G)$. By Deligne \cite{D0} Cor. 4.2.9 we
know that $\mathfrak{g}$ is semi-simple.  We then complexify
$\chi$ to obtain $\chi_{\Complex}: \mathfrak{g}_{\Complex}\to
\SymplecticLie(20,\Complex)$,
which is irreducible in the sense that after composition with the
natural representation
$$
\SymplecticLie(20,\Complex)\longrightarrow \mathfrak{gl}(20,\Complex)
$$
$\chi_{\Complex}$ is an irreducible representation of the
semi-simple complex Lie algebra $\mathfrak{g}_{\Complex}$.\\

\noindent
{\em Step 2.} In this step we classify all possible complex Lie algebra
monomorphism
$\chi_{\Complex}: \mathfrak{g}_{\Complex} \to \SymplecticLie(20,\Complex)$
where $\mathfrak{g}_{\Complex}$ is semi-simple and
$\chi_{\Complex}$ is irreducible in the sense described above. In order to
classify
$(\mathfrak{g}_{\Complex},\chi_{\Complex})$, we observe that it
suffices to consider all 20-dimensional irreducible
representations of complex semi-simple Lie algebras
$\mathfrak{g}_{\Complex}$. Actually, an irreducible representation
$\mathfrak{g}_{\Complex}\to \mathfrak{gl}(V_{\Complex})$ with
$\dim(V_{\Complex})>20$ admits no factorization
$$
\mathfrak{g}_{\Complex}\longrightarrow \SymplecticLie(20,\Complex)
\longrightarrow \mathfrak{gl}(V_{\Complex}).
$$
The reason is that, since $\mathfrak{g}_{\Complex}$ is mapped onto a proper
subspace of $\SymplecticLie(20,\Complex)$, the composition must decompose and
hence is reducible. We can list all such possibilities. Our method is
first to find all 20-dimensional representation of a semi-simple
Lie algebra, and then exclude those whose images do not lie in
$\SymplecticLie(20,\Complex)$. \smallskip

\noindent
{\em Case 1.} $\mathfrak{g}_{\Complex}$ has only one simple factor:
\begin{itemize}
  \item[(a)] $(A_1,[19])$,
  \item[(b)] $(A_5,[0,0,1,0,0])$,
  \item[(c)] $(C_2,[3,0])$.
\end{itemize}

\noindent
{\em Case 2.} $\mathfrak{g}_{\Complex}$ has two simple factors:
\begin{enumerate}
  \item[(a)] $(A_1\oplus C_2,[1]\otimes[2,0])$,
  \item[(b)] $(A_1\oplus C_2,[4]\otimes[1,0])$,
  \item[(c)] $(A_1\oplus D_5,[1]\otimes[1,0,0,0,0])$,
  \item[(d)] $(C_2\oplus C_2,[1,0]\otimes[0,1])$.
\end{enumerate}
A 20-dimensional representation has image in $\SymplecticLie(20,\Complex)$ if
and only if there exists an one-dimensional component in the
irreducible decomposition of the second wedge power. Here is an
example. The pair $(A_1\oplus A_1\oplus A_4,[1]\otimes
[1]\otimes[1,0,0,0])$ associates to the semi-simple Lie algebra
$A_{1}\oplus A_1\oplus A_4$ a 20-dimensional representation with
the highest weight $[1]\otimes [1]\otimes[1,0,0,0]$. One easily
checks that in the irreducible decomposition
\begin{eqnarray*}
  \bigwedge^2([1]\otimes [1]\otimes[1,0,0,0]) &\simeq& [0]\otimes [2]\otimes[2,0,0,0]\oplus [2]\otimes
[0]\otimes[2,0,0,0] \\
   &\oplus & [0]\otimes [0]\otimes[0,1,0,0] \oplus [2]\otimes
   [2]\otimes[0,1,0,0],
\end{eqnarray*}
there is no one dimensional component.\\

\noindent
{\em Step 3}.  All possible simple real groups of Hodge types are listed
in \S4 \cite{Sim}. Based on this and the classification given in
the last step we can now discuss them case by case by applying
the plethysm method. However, the following general result about
the $\Complex$-PVHS structures on a tensor product will simplify
our argument to a large extent. Since this result is of
interest in itself, we would like to include a proof in this paper.\\

Let $\bb{V}$ be an irreducible $\Complex$-PVHS over a
quasi-projective manifold $X\setminus S$ and with unipotent local
monodromy around $S$. Let
$$
\rho: \pi_{1}(X\setminus S)\to \GeneralLinear(V)
$$
be the corresponding representation of the fundamental group and
$G$ be the Zariski closure of $\rho$. Assume $G=G_{1}\times G_{2}$
with $G_{i}$ simple. Then  according to Schur's lemma $\bb{V}$ is
decomposed into
$$ \bb{V}\simeq \bb{V}_{1}\otimes \bb{V}_{2},$$
where $\bb{V}_{i}$ corresponds to a representation
$$
\rho_i: \pi_1(X\setminus S)\to G_i.
$$
\begin{proposition}\label{Tensorproduct}
The $\Complex$-PVHS on $\bb{V}$ factors into PVHS's on each factor
$\bb{V}_{i}$, i.e.\ each $\bb{V}_i$ admits a $\Complex$-PVHS
structure such that their tensor product on $\bb{V}_{i}$ coincides
with the $\Complex$-PVHS on $\bb{V}$.
\end{proposition}

\noindent
{\em Proof of Proposition 5.4:} We write $\dim V_i=n_i$ for
$i=1,2$ and we assume that $n_1\geq n_2$ without lose of
generality. We first need the following lemma.
\begin{lemma}\label{unipotent}
Each factor $\rho_i$ has  quasi-unipotent local monodromy around
$S$.
\proof By choosing a base point in $S$, the tensor product decomposition of
$\bb{V}$ gives the tensor product decomposition of the vector space
$V\simeq V_1\otimes V_2$ with group action, and since $G_i$ is simple,
$G_i\subset \SpecialLinear(V_i)$ for $i=1,2$. Now we apply $\bigwedge^{n_2}$
on the above isomorphism. Ex.\ 6.11(b) in \cite{FH} tells us that, for $V$
considered as a representation space of $\SpecialLinear(V_1)
\times \SpecialLinear(V_2)$, there exists an
irreducible component
$$
S^{n_2}(V_1)\subset \bigwedge^{n_2}(V).
$$
Since $\bb{V}$ is of unipotent local monodromy, each direct
component of $\bigwedge^{n_2}(V)$ is of unipotent local monodromy,
too. In particular, $S^{n_2}(V_1)$ is of unipotent local
monodromy. Let $T$ be one of local monodromy operators of
$\rho_1$, and $\lambda$ be one of eigenvalues of $T$. Then
clearly, $\lambda^{n_2}$ is one of eigenvalues of $T$ on
$S^{n_2}(V_1)$, hence is equal to one. This proves that $\rho_1$ is of
quasi-unipotent local monodromy. And by the unipotency of $\rho$,
$\rho_2$ is of quasi-unipotent local monodromy as well. This
completes the proof of Lemma \ref{unipotent}.\qed
\end{lemma}

Since $\rho_i: \pi_1(X\setminus S)\to G_i$ is a Zariski dense
representation into the simple algebraic group $G_i$ and with
quasi-unipotent local monodromy around $S$, by Jost-Zuo \cite{JZ}
there exists a pluri-harmonic metric on the flat bundle $\bb{V}_i$
with finite energy, which makes $\bb{V}_i$ into a Higgs bundle
$(E,\theta)_i$ over $X\setminus S$. Furthermore, T. Mochizuki
\cite{M} has analyzed the singularity of this harmonic metric in
detail and has shown that $(E,\theta)_i$ admits a logarithmic
extension $(\bar E,\bar\theta)_i$ over $X,$ i.e.\ $\bar E_i$ is an
extension of $E_i$, $\bar \theta_i$ is an extension of $\theta_i$
and such that
$$
\bar\theta: \bar E_i\to \bar E_i\otimes\Omega^1_X(\log S).
$$
Such a pluri-harmonic metric is called {\em tame}. In this case the
residue of $\bar\theta$ along $S$ is nilpotent. \\

From the proof of Lemma \ref{unipotent}, we know that, by applying
the Schur functor $\bigwedge^{n_2}$, one finds a direct factor of
$\bigwedge^{n_2}(\bb{V}_1\otimes\bb{V}_2)$ of the form
$$
S^{n_2}(\bb{V}_1)\otimes\det(\bb{V}_2),
$$
and  $S^{n_2}$ is non-trivial. Since  $G_2$ is simple, $\det(\bb{V}_2)$ is the trivial representation. \\

We consider $G_1$ as a simple algebraic subgroup of $\GeneralLinear(V_1).$
Since the Schur functor $S^{n_2}$ is non-trivial and $G_1$ is a
simple algebraic group, the representation
$$
S^{n_2}: G_1\to \GeneralLinear(S^{n_2}(V_1))
$$
is faithful. Since $\rho_1$ is Zariski dense in $G_1$,
$S^{n_2}(\rho_1)$ is irreducible. Since
$\bigwedge^{n_2}(\bb{V}_1\otimes\bb{V}_2)$ is semi-simple, there
exists a decomposition
$$
\bigwedge^{n_2}(\bb{V}_1\otimes\bb{V}_2)=\bigoplus_{i} S_i\otimes
W_i,
$$
where $S_i$ are irreducible and $W_i$ are trivial.  By Deligne's
Prop.\ 1.13 in \cite{D1}, there exists uniquely $\Complex$-PVHS on
$S_i$ and $\Complex$-HS on $W_i$ such that the direct sum of the
tensor products of them coincides with the $\Complex$-PVHS on
$\bigwedge^{n_2}(\bb{V}_1\otimes\bb{V}_2)$. So, in particular,
there exists a $\Complex$-PVHS on $S^{n_2}(\bb{V}_1).$ \\

By the uniqueness of such pluri-harmonic metric, $S^{n_2}(\bar
E,\bar\theta)_1$ coincides with the $\Complex$-PVHS on
$S^{n_2}(\bb {V}_1)$. Hence $S^{n_2}(\bar E,\bar\theta)_1$ is a
fixed point of the $\Complex^\times$-action.
The representation $G_1\to \GeneralLinear(S^{n_2}(V_1))$ induces a morphism
$$
\phi_{S^{n_2}}: M(\pi_1(X\setminus S), G_1)^{s.s}\to
M(\pi_1(X\setminus S),\GeneralLinear(S^{n_2}(V_1)))^{s.s}
$$
between the corresponding moduli spaces of semi-simple
representations. By Simpson's Cor.\ 9.16 in\cite{Sim1}, $\phi_{S^{n_2}}$
is finite. \\

If $S=\emptyset,$ then $\Complex^\times$ acts on both moduli spaces
continuously via Hermitian-Yang-Mills metrics on poly-stable Higgs
bundles  $(E, t\theta)$. And this action is compatible with
$\phi_{S^{n_2}}.$  Since $S^{n_2}(\rho_1)$ is a fixed point of
$\Complex^\times$-action, $\rho_1$ is a fixed point of
$\Complex^\times$-action. This just means that $(E,\theta)_1$ is a
fixed point of $\Complex^\times$-action. Hence $(E,\theta)_1$ is a
$\Complex$-PVHS on $\bb{V}_1.$ In general $S\not=\emptyset.$ We
take a curve $C\setminus S\subset X\setminus S$, which is a
complete intersection of ample hypersurfaces. Taking the
restrictions
$$
\rho_1|_{C\setminus S}\in M(\pi_1(C\setminus S), G_1)^{s.s},
$$
we have
$$
S^{n_2}(\rho_1)|_{C\setminus S}\in M(\pi_1(C\setminus S),
\GeneralLinear(S^{n_2}(V_1)))^{s.s}.
$$
We consider the map
$$
\phi_{S^{n_2}}: M(\pi_1(C\setminus S), G_1)^{s.s}\to
M(\pi_1(C\setminus S), \GeneralLinear(S^{n_2}(V_1)))^{s.s}.
$$
By Simpson's main theorem in \cite{Sim0}, there exist
Hermitian-Yang-Mills metrics on poly-stable Higgs bundles on $C$
with logarithmic pole of the Higgs field on $S$. And the
$\Complex^\times$-action can be defined on both spaces of semi-simple
representations on $C\setminus S$ via Hermitian-Yang-Mills metric on
$(\bar E,t\bar\theta).$ Applying the same argument as above to the
compact case, we show that the pulled back Higgs bundle $(\bar
E,\bar \theta)_1$ to $C\setminus S$ is a fixed point of the
$\Complex^\times$-action.  If we choose $C\setminus S$ sufficiently
ample, then $(\bar E,\bar\theta)_1$ is also a fixed point of the
$\Complex^\times$-action. (The isomorphism $(\bar
E,\bar\theta)_1|_C\simeq(\bar E,t\bar\theta)_1|_C$ extends to an
isomorphism $(\bar E,\bar\theta)_1\simeq(\bar E,t\bar\theta)_1$ if
$C$ is sufficiently ample.) Again by Simpson, $(\bar E,\bar\theta)_1$ is a $\Complex$-PVHS on $\bb{V}_1$. \\

Similarly, we also show that $\bb{V}_2$ admits a $\Complex$-PVHS.
The tensor product of $\Complex$-PVHS on $\bb{V}_1$ and on
$\bb{V}_2$ is a $\Complex$-PVHS on $\bb{V}_1\otimes\bb{V}_2.$ By
Deligne's uniqueness theorem on $\Complex$-PVHS on irreducible
local systems, this tensor product coincides with the original
$\Complex$-PVHS on $\bb{V}_1\otimes\bb{V}_2.$  Proposition
\ref{Tensorproduct} is completed. \qed \\

Now we start with the analysis of case 2. By the above
proposition, we know that in this case we have an isomorphism
$$
(E,\theta)\simeq (E_1,\theta_1)\otimes (E_2,\theta_2),
$$
where each $(E_i,\theta_i)$ is a system of Hodge bundles. Because
the Hodge numbers of $E$ are 1,9,9,1, it is not difficult to see
that, up to permutation of factors, the Hodge numbers of
$(E_1,\theta_1)$ are 1,1, and those of $(E_2,\theta_2)$ are 1,8,1.
Since $(E_1,\theta_1)$ comes from a $\Reals$-PVHS,
$(E_2,\theta_2)$ also comes from a $\Reals$-PVHS. This implies
that, if $\mathfrak{g}=\mathfrak{g}_1\oplus \mathfrak{g}_{2}$,
then up to permutation of factors,
$\mathfrak{g}_1=\SpecialUnitaryLie(1,1)$ and
$\mathfrak{g}_2\subset \SpecialOrthogonalLie(2,8)$.\\

From the list in \S4 \cite{Sim}, we know that the possible real
forms of $\SymplecticLie(4,\Complex)$ are
$\SymplecticLie(1,1)$ and $\SymplecticLie(4,\Reals)$,
and those of $\SpecialOrthogonalLie(10,\Complex)$ are
$\SpecialOrthogonalLie(2,8)$ and $\SpecialOrthogonalLie(4,6)$,
respectively. It is straightforward to check that the only 
possibility of case 2 is
\begin{itemize}
\item case (2c) \quad $(\SpecialUnitaryLie(1,1)\oplus
\SpecialOrthogonalLie(2,8), \id\otimes\id)$.
\end{itemize}
Obviously, case (1a) is impossible since it is of weight 19. For
those real forms of Hermitian types in case 1 we can again compare
the Hodge numbers. The only possibilities are
\begin{itemize}
    \item case (1b) \quad $(\SpecialUnitaryLie(3,3),\bigwedge^3)$;
    \item case (1c) \quad $(\SymplecticLie(1,1),S^3)$.
\end{itemize}
Note that case (1c) is of non-Hermitian type. \\

\noindent
{\em Step 4.} We have already excluded case (1b) using the method of
plethysm. In the last step, we apply the method further in order to exclude
the left two cases. The argument for case (2c) is similar. We give
the analogous statement of Proposition \ref{PlethysmProp} as
follows:
$$
S^{2}(E,\theta)=(E_1,\theta_1)\oplus (E_2,\theta_2)\oplus
(E_3,\theta_3)
$$
where
\begin{eqnarray*}
E_{1} &=&\phantom{E_{1}^{6,0}\oplus E_{1}^{5,1}\oplus}\
E_{1}^{4,2}\oplus
E_{1}^{3,3}\oplus E_{1}^{2,4} \phantom{E_{1}^{1,5}\oplus E_{1}^{0,6}}  \\
E_{2} &=&\phantom{E_{1}^{6,0}\oplus E_{1}^{5,1}\oplus}\
E_{2}^{4,2}\oplus
E_{2}^{3,3}\oplus E_{2}^{2,4} \phantom{E_{1}^{1,5}\oplus E_{1}^{0,6}}  \\
E_{3} &=& E_{3}^{6,0}\oplus E_{3}^{5,1}\oplus E_{3}^{4,2} \oplus
E_{3}^{3,3} \oplus E_{3}^{2,4}\oplus E_{3}^{1,5}\oplus
E_{3}^{0,6}.
\end{eqnarray*}
The Hodge numbers of $E_{1}$ are $0, 0, 1, 1, 1, 0, 0$,
respectively, those of $E_{2}$ are $0, 0, 8$, $29, 8, 0, 0$ and those of $E_3$
are $1, 9, 45, 52, 45, 9, 1$. So by the computational result in Prop.\
\ref{PlethysmCompThm}, we see that case (2c) is impossible. For
case (1c), the corresponding result is the following:
$$
S^{2}(E,\theta)=(E_1,\theta_1)\oplus (E_2,\theta_2)\oplus
(E_3,\theta_3)\oplus (E_4,\theta_4)
$$
where the respective dimensions of $E_i$ are $10,35,81,84$. But we
are unable to obtain further information on the Hodge numbers of
$E_i$, because $\SymplecticGroup(1,1)$ is of non-Hermitian type.
But fortunately we can still get a contradiction by the actual
computation. The argument works as follows. The first Hodge bundle of
$S^2(E,\theta)$ is of dimension $1$, it must lie in one of $(E_i,\theta_i)$, and
hence the rank of the Higgs subsheaf generated by the first Hodge
bundle will not exceed $84$. Over the point used in Prop.\
\ref{PlethysmCompThm}, the rank of the stalk of the Higgs subsheaf
generated by the first Hodge bundle is not less than
$$
1+9+45+78=133.
$$
This gives the desired contradiction for case (1c). The proof is
complete.\qed


\begin{thebibliography}{alpha}
\addcontentsline{toc}{chapter}{Bibliography}

\bibitem{ACT}
Allcock,D; Carlson,J.; Toledo,D; {\it The complex hyperbolic
geometry of the moduli space of cubic surfaces}.  J. Algebraic
Geometry.  11 (2002),  no. 4, 659-724.

\bibitem{Periods}
Carlson,J; M\"uller-Stach,S; Peters,C; {\it Period Mappings and
Period Domains}, Cambridge Studies in Advanced Mathematics,
Cambridge University Press 2003.

\bibitem{D0}
Deligne,P; {\it Th\'eorie de Hodge II}, Publ.Math.I.H.E.S.
40(1971), 5-57.

\bibitem{D1}
Deligne,P; {\it Un th\'{e}or\`{e}me de finitude pour la
monodromie}, Discrete Groups in Geometry and Analysis, Birkhauser
1987, 1-19.

\bibitem{DKV}
Dolgachev,I.; van Geemen,B.; Kond\=o,S; {\it A complex ball
uniformization of the moduli space of cubic surfaces via periods
of K3 surfaces}.  J. Reine Angew. Math.  588  (2005), 99-148.

\bibitem{FH}
Fulton,W; Harris,J; {\it Representation theory, A first course},
GTM 129.

\bibitem{G}
Griffiths,P; {\it Periods of integrals on algebraic manifolds II.
Local study of the period mapping}, American J. Math. (90), 1968
805-865.


\bibitem{JZ}
Jost,J; Zuo,K; {\it Harmonic maps and ${\rm
Sl}(r,C)$-representations of fundamental groups of quasiprojective
manifolds}.  J. Algebraic Geometry.  5 (1996),  no. 1, 77-106.

\bibitem{M}
Mochizuki,T; {\it Asymptotic behaviour of tame nilpotent harmonic
bundles with trivial parabolic structure}.  J. Differential
Geometry. 62  (2002),  no. 3, 351-559.

\bibitem{Mok}
Mok, N;{\it Uniqueness theorems of Hermitian metrics of
seminegative curvature on quotients of bounded symmetric domains}.
Annals of Mathematics. 125(1987), No.1, 105-152.

\bibitem{MSY}
Matsumoto,K; Sasaki,T; Yoshida,M; {\it The monodromy of the period
map of a 4-parameter family of K3 surfaces and the hypergeometric
function of type (3,6)}, International Journal of Mathematics, Vol
3, No.1, 1-164, 1992.


\bibitem{MVZ}
Moeller,M; Viehweg,E; Zuo,K; {\it Stability of Hodge bundles and a
numerical characterization of Shimura varieties}, arXiv
AG/07063462.


\bibitem{Nagel}
Nagel,J; {\it The image of the Abel-Jacobi map for complete
intersections}, Ph.D Thesis, Rijksuniversiteit Leiden 1997.

\bibitem{S}
Sheng,M; {\it On the geometric realizations of Hermitian symmetric
domains}, Ph.D Thesis, The Chinese University of Hong Kong 2005.

\bibitem{Sim0}
Simpson,C; {\it Harmonic bundles on noncompact curves}. J. Amer.
Math. Soc. 3 (1990), no. 3, 713-770.

\bibitem{Sim}
Simpson,C; {\it Higgs bundles and local systems}, Publ. Math.
I.H.E.S., No.75, 5-95, 1992.

\bibitem{Sim1}
Simpson,C; {\it Moduli of representations of the fundamental group
of a smooth projective varieites II}, Publ. Math. I.H.E.S., No.80,
5-79, 1994.

\bibitem{SYY}
Sasaki,T; Yamaguchi,K; Yoshida,M; {\it On the ridity of
differential systems modelled on Hermitian symmetric spaces and
disproofs of a conjecture concerning modular interpretations of
confinguration spaces}. Advanced Studies in Pure Mathematics 25,
318-354, 1997.

\bibitem{SZ}
Sheng,M; Zuo,K; {\it Calabi-Yau like PVHS and characteristic
subvariety over bounded symmetric domains}, arXiv: AG/07053779.

\bibitem{VZ1}
Viehweg,E; Zuo,K; {\it A characterization of certain Shimura
curves in the moduli stack of abelian varieties}, J. Differential
Geometry, 66, No.2, 233-287, 2004.

\bibitem{VZ2}
Viehweg,E; Zuo,K; {\it Arakelov inequalities and the
uniformization of certain rigid Shimura varieties}, To be appeared
in J. Differential Geometry.



\end{thebibliography}
\end{document}